\documentclass{amsart}

\usepackage{amssymb}
\usepackage{xcolor}
\definecolor{darkgreen}{rgb}{0,0.45,0} 
\definecolor{darkishgreen}{rgb}{0,0.65,0} 
\definecolor{darkred}{rgb}{0.75,0,0}
\definecolor{darkblue}{rgb}{0,0,0.6} 
\usepackage[pdfborder=0,pagebackref,colorlinks,citecolor=darkgreen,linkcolor=darkgreen,urlcolor=darkblue]{hyperref}
\renewcommand*{\backref}[1]{}
\renewcommand*{\backrefalt}[4]{({%
    \ifcase #1 Not cited.%
          \or on p.~#2%
          \else on pp.~#2%
    \fi%
    })}
\usepackage{mathtools}
\usepackage{comment}
\usepackage{enumitem}
\usepackage{tikz}
\usepackage{tikz-cd}

\newcommand{\id} { \mathrm{id}}
\newcommand{\ob} { \mathrm{ob}}

\newcommand{\dom}{\mathrm{dom}}
\newcommand{\cod}{\mathrm{cod}}

\newcommand{\Hom}{\mathrm{Hom}}

\newcommand{\To}{\Rightarrow}

\DeclareSymbolFont{bbold}{U}{bbold}{m}{n}
\DeclareSymbolFontAlphabet{\mathbbb}{bbold}
\newcommand{\catone}{\ensuremath{\mathbbb{1}}}
\newcommand{\cattwo}{\ensuremath{\mathbbb{2}}}

\newcommand{\cat}[1]{\textup{\textsf{#1}}}
\newcommand{\Cat}{\mathbb{C}\cat{at}}
\newcommand{\Bicat}{\mathbb{B}\cat{icat}}
\newcommand{\sSet}{\cat{s}\mathbb{S}\cat{et}}

\newcommand{\twoCat}{2\text{-}\mathbb{C}\cat{at}}

\newcommand{\Icon}{\mathbb{I}\cat{con}}

\newcommand{\cA}{{\mathcal{A}}}
\newcommand{\cB}{{\mathcal{B}}}
\newcommand{\cC}{{\mathcal{C}}}

\newcommand{\cJ}{{\mathcal{J}}}
\newcommand{\cX}{{\mathcal{X}}}

\newtheorem{thm}{Theorem}[section]

\newtheorem{prop}[thm]{Proposition}
\newtheorem{lem}[thm]{Lemma}

\theoremstyle{definition}
\newtheorem{defn}[thm]{Definition}

\newtheorem{obs}[thm]{Observation}
\newtheorem{rmk}[thm]{Remark}

\newtheorem{ex}[thm]{Example}

\newtheorem{clm*}{Claim}

\theoremstyle{remark}

\makeatletter
\let\c@equation\c@thm
\makeatother
\numberwithin{equation}{section}

\title{On \texorpdfstring{$\infty$}{infinity}-cosmoi of bicategories}

\author{Emily Riehl}
\email{eriehl@jhu.edu}
\address{Department of Mathematics, Johns Hopkins University, Baltimore, MD 21218}

\author{Mira Wattal}
\email{wattal@bc.edu}
\address{Mathematics Department, Boston College, Chestnut Hill, MA 02467}

\thanks{This work was supported by the NSF via DMS-1652600, by the ARO under MURI Grant W911NF-20-1-0082, and by the Johns Hopkins President's Frontier Award program. The authors also wish to thank John Bourke, Alexander Campbell, and Dominic Verity for enlightening conversations related to this project, and further thank Alexander Campbell for comments and corrections on the first draft.}

\begin{document}

\maketitle

\begin{abstract} An $\infty$-cosmos is a setting in which to develop the formal category theory of $(\infty,1)$-categories. In this paper, we explore a few atypical examples of $\infty$-cosmoi whose objects are 2-categories or bicategories rather than $(\infty,1)$-categories and compare the formal category theory that is so-encoded with classical 2-category and bicategory theory. We hope this work will inform future explorations into formal $(\infty,2)$-category theory.
\end{abstract}

\setcounter{tocdepth}{1}
\tableofcontents

\section{Introduction}

In \cite{Street:1974ec}, Street introduced the notion of an \emph{elementary cosmos}, axiomatizing the 2-categorical ``universe'' in which categories live as objects. Inspired by this work, Riehl and Verity developed a corresponding notion of an \emph{$\infty$-cosmos}, axiomatizing the $(\infty,2)$-categorical ``universe'' in which $(\infty,1)$-categories live as objects; see \cite{RV-elements} summarizing a previous series of papers. The objects of an $\infty$-cosmos are called ``$\infty$-categories'' --- turning the nickname Lurie popularized in his book \cite{Lurie:2009ht} into a technical term --- and the morphisms between them are called ``$\infty$-functors.'' Any $\infty$-cosmos has a quotient 2-category, called the \emph{homotopy 2-category}, whose objects and 1-cells are the $\infty$-categories and $\infty$-functors in the $\infty$-cosmos and whose 2-cells then define ``$\infty$-natural transformations'' between them.

 From the axioms of an $\infty$-cosmos, Riehl and Verity define the notions of equivalences and adjunctions between its objects --- that is, equivalences and adjunctions between $\infty$-categories --- and limits and colimits of diagrams valued inside an object of an $\infty$-cosmos --- that is, limits and colimits inside $\infty$-categories --- and prove that these recover the standard definitions from the $(\infty,1)$-categorical literature. They then provide new formal proofs of expected categorical theorems: for instance that right adjoints between $\infty$-categories preserve limits inside those $\infty$-categories. The aim of this program is to provide a ``model-independent'' foundation of $(\infty,1)$-category theory.
 
 The axioms of an $\infty$-cosmos are deliberately very minimal. Consequently, in addition to the expected examples --- $\infty$-cosmoi whose objects are $(\infty,1)$-categories in some model --- there are a variety of ``exotic'' examples of $\infty$-cosmoi. In particular, several models of $(\infty,2)$-categories or even $(\infty,n)$-categories give rise to $\infty$-cosmoi  \cite[\S E.3]{RV-elements}, including Ara's $n$-quasi-categories \cite{Ara:2014hq}, Rezk's $\Theta_n$-spaces \cite{Rezk:2010cp}, Barwick's $n$-fold complete Segal spaces \cite{Barwick:2005in}, and Verity's $n$-complicial sets \cite{Verity:2008wcI}. For each of these models and for each $0 \leq n < \infty$ (including $n=\infty$ in the final case), there is an $\infty$-cosmos in which the ``$\infty$-categories'' are the $(\infty,n)$-categories in that particular model. This suggests the tantalizing possibility that it might be possible to develop $(\infty,2)$-category theory or $(\infty,n)$-category theory ``model-independently'' by adapting $\infty$-cosmological methods.
 
Of course, all of the theorems proven about arbitrary $\infty$-cosmoi apply to every particular example. Thus, for instance, \cite[2.4.2]{RV-elements} proves that right adjoints between 2-quasi-categories preserve limits in 2-quasi-categories. What is not clear, however, is what exactly the $\infty$-cosmological definitions of adjunctions between $\infty$-categories or limits inside $\infty$-categories compile out to in the 2-quasi-categories model.

This motivates the present paper. In addition to the $\infty$-cosmoi of $(\infty,2)$-cat\-e\-gories mentioned above, there are --- somewhat curiously --- various $\infty$-cosmoi whose ``$\infty$-categories'' are 2-categories or bicategories and whose ``$\infty$-functors'' and ``$\infty$-nat\-u\-ral transformations'' define some variety of functor and natural transformation between them. In \S\ref{sec:cosmoi}, we develop  two particularly fertile examples: the $\infty$-cosmos $\twoCat$ of 2-categories, 2-functors, and 2-natural transformations and the $\infty$-cosmos $\Icon$ of bicategories, normal pseudofunctors, and icons.

In \S\ref{sec:formal}, we review the basic formal theory of $\infty$-categories in an $\infty$-cosmos, recalling the $\infty$-cosmological definitions of equivalences between $\infty$-categories, adjunctions between $\infty$-categories, terminal objects in an $\infty$-category, and limits in an $\infty$-category. In \S\ref{sec:2-cat-cosmos} and \ref{sec:icon-cosmos}, we explore each of the $\infty$-cosmoi $\twoCat$ and $\Icon$ in turn and describe the varieties of formal 2-category theory and bicategory theory they encode. We conclude in \S\ref{sec:conjectures} with a proposal for further study.
 
\section{\texorpdfstring{$\infty$}{Infinity}-cosmoi of 2-categories and bicategories}\label{sec:cosmoi}

The informal definition of an $\infty$-\emph{cosmos} states that it is an  ``$(\infty,2)$-category with $(\infty,2)$-categorical limits,'' under a particularly strict interpretation of these terms. Firstly, an $\infty$-cosmos is a simplicially-enriched category whose hom-spaces are quasi-categories, these being a model of $(\infty,1)$-categories. An $\infty$-cosmos must then possess certain simplicially-enriched limits including a terminal object, products, and simplicial cotensors. More generally, an $\infty$-cosmos possesses all simplicially-enriched limits whose weights are ``flexible.''\footnote{As noted in \cite[6.2.7]{RV-elements}, the term ``flexible weighted limits'' is a bit of a misnomer, as they are more precisely analogous to the PIE-limits in 2-category theory than the flexible ones, since splittings of idempotents are not required.} These flexible weighted limits can be built inductively as countable inverse limits of towers of pullbacks of products of a special class of maps in an $\infty$-cosmos called ``isofibrations.'' The class of isofibrations is then required to satisfy various closure properties; see  \cite[1.2.1]{RV-elements} for full details.

As noted in \cite[E.1.6]{RV-elements}, any 2-category with sufficient limits defines an $\infty$-cos\-mos, whose homotopy 2-category of ``$\infty$-cat\-e\-gories,'' ``$\infty$-functors,'' and ``$\infty$-natural transformations'' is isomorphic to the original 2-category.  Morally, a 2-category defines an $\infty$-cosmos just when that 2-category admits PIE-limits: products, inserters, and equifiers. However, as discussed in \cite[E.1.6]{RV-elements}, slightly more is required. An $\infty$-cosmos is required to come with a specified class of 1-cells, the above-mentioned ``isofibrations.'' This class must satisfy the closure properties enumerated in \cite[1.2.1(ii)]{RV-elements}. In particular, the 2-category is required to admit 2-pullbacks and inverse 2-limits of towers of these maps, not merely the corresponding bilimits that can be constructed as PIE-limits. 

An additional requirement of the class of isofibrations in an $\infty$-cosmos is that they also define isofibrations in the homotopy 2-category, in a sense we now recall, although in general the isofibrations in an $\infty$-cosmos may be a smaller class than just those maps that define isofibrations in the homotopy 2-category.

\begin{defn}\label{defn:isofibration} A 1-cell $F \colon \cA \to \cB$ in a 2-category defines an \textbf{isofibration} just when it possesses the lifting property for invertible 2-cells depicted below:
\begin{equation}\label{eq:isofibration}
\begin{tikzcd} \cX \arrow[r, bend left, "A"] \arrow[dr, bend right, "B"'] & \cA \arrow[d, "F"] & \cX \arrow[r, bend left, "A"] \arrow[r, "L"', dashed, bend right]  \arrow[r, phantom, "\scriptstyle\exists\cong\Downarrow\alpha"] & \cA \arrow[d, "F"] \\ \arrow[ur, phantom, "\scriptstyle\cong\Downarrow\beta" near end] & \cB \arrow[ur, phantom, "="] & & {~}\cB .
\end{tikzcd}
\end{equation}
Equivalently, the isofibrations in a 2-category can be characterized as ``representably-defined isofibrations.'' A 1-cell  $F \colon \cA \to \cB$ is an isofibration just when, for any object $\cX$, the induced functor between the hom-categories  is an isofibration of categories, with the right lifting property 
\[
\begin{tikzcd} \catone \arrow[r, "A"] \arrow[d, hook] & \Hom(\cX,\cA) \arrow[d, "F"] \\ \mathbb{I} \arrow[r, "\beta"'] \arrow[ur, dashed, "\alpha"] & \Hom(\cX,\cB)
\end{tikzcd}
\]
against the inclusion of either endpoint of the free-living isomorphism $\mathbb{I}$.
\end{defn}

For example, the 2-category $\Cat$ of ordinary categories, functors, and natural transformations is an $\infty$-cosmos in which the isofibrations are the isofibrations between categories, those functors that have the right lifting property with respect to $\catone\hookrightarrow \mathbb{I}$ \cite[1.2.11]{RV-elements}.

\subsection{2-categories, 2-functors, and 2-natural transformations}

Our first example might be thought of as an $\infty$-cosmos for strict 2-category theory:

\begin{prop}\label{prop:2cat-cosmos} The 2-category $\twoCat$ of 2-categories, 2-functors, and 2-natural transformations defines an $\infty$-cosmos whose isofibrations are the isofibrations.
\end{prop}

Note a 2-functor is an isofibration in $\twoCat$ just when its underlying functor is an isofibration in $\Cat$.

\begin{proof}
This 2-category, like any 2-category with  limits and colimits, admits a ``trivial model structure'' due to Lack. In this model structure, the fibrations are the isofibrations and the weak equivalences are the equivalences, as defined internally to the 2-category of 2-categories \cite[3.3]{Lack-2monad} (see Definitions  \ref{defn:isofibration} and \ref{defn:equivalence}). All objects are fibrant and cofibrant. Moreover, this model structure is enriched over the folk model structure on categories. Thus we can apply \cite[E.1.4]{RV-elements} using the adjunction
\[
\begin{tikzcd} 
\Cat \arrow[r, bend right, "N"'] \arrow[r, phantom, "\bot"] & \sSet \arrow[l, bend right, "\text{ho}"']
\end{tikzcd}
\]
to obtain an $\infty$-cosmos of 2-categories, in which the isofibrations are the isofibrations.
\end{proof}

\begin{rmk} 
The $\infty$-cosmos $\twoCat$ is cartesian closed as an $\infty$-cosmos. It satisfies the additional axiom of \cite[1.2.23]{RV-elements}, which in this context amounts to the following observations:
\begin{enumerate}
\item The 2-category $\cB^\cA$ of 2-functors, 2-natural transformations, and modifications defines a right 2-adjoint to the cartesian product.
\item Exponentiation preserves isofibrations.
\end{enumerate}
\end{rmk}

\subsection{Bicategories, normal pseudofunctors, and icons}

The strict 2-functors are just one variety of homomorphisms between 2-categories or bicategories. Arguably more natural are the \emph{pseudofunctors}, which preserve composition and identities up to coherent isomorphism, or the \emph{normal pseudofunctors}, which preserve identities strictly and composition up to coherent isomorphism. Lack proves that either class of functors defines the 1-cells of a 2-category whose 2-cells are the \emph{icons}, the identity-component oplax natural transformations \cite[3.1]{Lack-icons}.

\begin{defn}\label{defn:icon} For normal pseudofunctors $F,G \colon \cA \to \cB$ satisfying $Fx = Gx$ for all objects $x \in \cA$, an \textbf{icon} $\alpha \colon F \To G$ is given by the data of a natural transformation
\[
\begin{tikzcd} \cA(x,y) \arrow[r, bend left, "F"] \arrow[r, bend right, "G"'] \arrow[r, phantom, "\scriptstyle\Downarrow\alpha_{x,y}"] & \cB(Fx,Fy)
\end{tikzcd}
\]
for all $x,y \in \cA$ satisfying two axioms:
\begin{enumerate}
\item The component of $\alpha_{x,x}$ at $\id_x \in \cA(x,x)$ is the identity 2-cell on $F(\id_x) = \id_{Fx} = \id_{Gx}= G(\id_x)$.
\item For all $x,y,z \in \cA$, 
\[
\begin{tikzcd}[row sep=large, column sep=1.5em]
\cA(y,z) \times \cA(x,y) \arrow[r, "\circ"] \arrow[d, bend right=60, "F \times F"'] \arrow[d, bend left=60, "G \times G"] \arrow[dr, phantom, "\scriptstyle\qquad\cong" pos=.45, "\scriptstyle\qquad\Rightarrow" pos=.6] \arrow[d, phantom, "\scriptstyle\alpha_{y,z} \times \alpha_{x,y}" pos=.4, "\scriptstyle\Rightarrow " pos=.6] & \cA(x,z) \arrow[d, bend left, "G"] \arrow[dr, phantom, "="] & [-15pt] \cA(y,z) \times \cA(x,y) \arrow[d, bend right, "F \times F"'] \arrow[r, "\circ"] \arrow[dr, phantom, "\scriptstyle\cong" pos=.45, "\scriptstyle\Rightarrow" pos=.6] & \cA(x,z) \arrow[d, bend right, "F"'] \arrow[d, bend left, "G"] \arrow[d, phantom, "\scriptstyle\alpha_{x,z}" pos=.4, "\scriptstyle\Rightarrow" pos=.6] \\ 
\cB(Fy,Fz) \times \cB(Fx,Fy) \arrow[r, "\circ"'] & \cB(Fx,Fz) & \cB(Fy,Fz) \times \cB(Fx,Fy) \arrow[r, "\circ"'] & \cB(Fx,Fz)
\end{tikzcd}
\]
where the unlabelled isomorphisms are the coherences of the normal pseudofunctors $F$ and $G$.
\end{enumerate}
\end{defn}

There are two independent choices we might make in defining a 2-category with PIE-limits whose 2-cells are icons. We can take the objects to be the 2-categories or the bicategories, and we can take the 1-cells to be the pseudofunctors or the normal pseudofunctors. In fact, all four possibilities are biequivalent \cite{LackPaoli-2Nerves}. However, while we can prove that either 2-category of bicategories defines an $\infty$-cosmos, the 2-categories of 2-categories appear not to, as discussed in Remark \ref{rmk:2cat-icon}.

\begin{prop} 
The 2-category $\Icon$ of bicategories, normal pseudofunctors, and icons defines an $\infty$-cosmos, in which the isofibrations are the isofibrations.
\end{prop}

Note a normal pseudofunctor is an isofibration in $\Icon$ just when all of the induced functors between hom-categories are isofibrations in $\Cat$.

\begin{proof}
The 2-category of bicategories, normal pseudofunctors, and icons is the 2-category of strict algebras, pseudo morphisms, and algebra 2-cells for a flexible 2-monad on the 2-category of reflexive $\Cat$-graphs. As such it admits all flexible limits. By \cite[A.1]{Bourke-accessible}, this 2-category then also admits 2-pullbacks of normal isofibrations, and by a similar argument, limits of towers of normal isofibrations. 

A normal isofibration is an isofibration admitting chosen lifts that are natural in $\cX$ and defined in such a way that the lifts of identity 2-cells are identities. In $\Icon$, all isofibrations are normal isofibrations, since all isofibrations in $\Cat$ are normal isofibrations. Thus, $\Icon$ defines an $\infty$-cosmos in which the isofibrations are the isofibrations.
\end{proof}

The 2-category of bicategories, pseudofunctors, and icons also defines an $\infty$-cos\-mos for a similar reason; in this case the flexible 2-monad is on the category of $\Cat$-graphs. We prefer to work with normal pseudofunctors as their data is easier to enumerate. For instance, a normal pseudofunctor $b \colon 1 \to \cB$ is just given by the data of an object $b \in \cB$, while a pseudofunctor $b \colon 1 \to \cB$ is given by an object $b \in \cB$, an endoequivalence $i \colon b \to b$, and an invertible 2-cell $\beta \colon \id_b \cong i$ satisfying a coherence condition.

\begin{rmk}\label{rmk:2cat-icon}
The full subcategory of 2-categories, normal pseudofunctors, and icons is nearly an $\infty$-cosmos. One can directly verify that this 2-category has PIE-limits inherited from $\Icon$. However, the 2-pullback of a normal isofibration is a bicategory, rather than a 2-category --- unless that normal isofibration happens to be \emph{discrete}, meaning that each invertible 2-cell has a unique lift. The fact that the 2-category of 2-categories, normal pseudofunctors, and icons admits 2-pullbacks of discrete isofibrations  suggests an alternate potential $\infty$-cosmos structure where the isofibrations are taken to be the discrete isofibrations. Alas, when a 2-category $C$ admits non-identity invertible 2-cells, the unique functor $C \to \catone$ is not a discrete isofibration. So this class violates one of the closure axioms required of an $\infty$-cosmos in \cite[1.2.1(ii)]{RV-elements}.
\end{rmk}

\subsection{Non-examples}

Many examples of $\infty$-cosmoi arise as full subcategories of fibrant objects in a suitably-enriched model category with all objects cofibrant in which case the isofibrations taken to be the fibrations between fibrant objects. This was the strategy employed in the proof of Proposition \ref{prop:2cat-cosmos}.  Such $\infty$-cosmoi are well-behaved in the sense that the equivalences between $\infty$-categories, as captured by Definition \ref{defn:equivalence}, are precisely the weak equivalences in the model structure \cite[1.2.13, 1.4.7]{RV-elements}. This suggests that we might productively consider model categories of 2-categories or bicategories.

In \cite{Lack-2cat, Lack-bicat}, Lack describes a pair of Quillen equivalent model structures on the categories of 2-categories and 2-functors and bicategories and strict functors, respectively. These are attractive settings within which to consider the homotopy theory of 2-categories or bicategories. In each case, all objects are fibrant, and the weak equivalences are the biequivalences that respectively define 2-functors or pseudofunctors. Unfortunately, in order to extract an interesting $\infty$-cosmos structure from either of these model categories we must find a compatible enrichment. By \cite[\S 7]{Lack-2cat}, the model structure for 2-categories is not cartesian monoidal, though it does define a monoidal model category with respect to the Gray tensor product. We are not aware of a strong monoidal Quillen adjunction relating this monoidal model category to the Joyal model structure on simplicial sets, and thus are unable to convert this into the required quasi-categorical enrichment by applying \cite[E.1.3]{RV-elements}. The model structure for bicategories fails to be cartesian monoidal for similar reasons.

\section{Formal category theory in an \texorpdfstring{$\infty$}{infinity}-cosmos}\label{sec:formal}

For the reader's convenience, we recall the following definitions from \cite[Chapter 2]{RV-elements}, which develop the basic formal category theory of ``$\infty$-categories,'' the objects in an $\infty$-cosmos. The context for the following definitions is the homotopy 2-category of $\infty$-categories, $\infty$-functors, and $\infty$-natural transformations. Here it suffices to work in any 2-category that has a 2-terminal object $1$ and is either cartesian closed or is cotensored over 1-categories.\footnote{In the examples we consider, these cotensors are defined strictly, though a weaker notion of cotensor would also suffice; see  \cite[\S 3.2]{RV-elements}.} A useful exercise for the reader is to verify that each of the definitions appearing below specializes to the standard categorical notions in the case where the $\infty$-cosmos is the 2-category $\Cat$, of ordinary categories, functors, and natural transformations.

Equivalences and adjunctions between $\infty$-categories are defined using the standard 2-categorical definitions.

\begin{defn}\label{defn:equivalence} An \textbf{equivalence}  is comprised of:
\begin{itemize}
    \item a pair of $\infty$-categories $\cA$ and $\cB$,
    \item a pair of $\infty$-functors $F\colon \cB \to \cA$ and $G \colon \cA \to \cB$,
    \item and a pair of invertible $\infty$-natural transformations $\beta \colon \id_\cB \cong GF$ and $\alpha \colon \id_\cA \cong FG$.
\end{itemize}
\end{defn}

\begin{defn}\label{defn:adjunction}
An \textbf{adjunction} is comprised of:
\begin{itemize}
\item a pair of $\infty$-categories $\cA$ and $\cB$,
\item a pair of $\infty$-functors $U \colon \cA \rightarrow \cB$ and $F \colon \cB \rightarrow \cA$,
\item and a pair of $\infty$-natural transformations $\eta \colon \id_\cB \Rightarrow UF$ and $\epsilon \colon FU \Rightarrow \id_\cA$, called the \textbf{unit} and \textbf{counit} respectively,
\end{itemize}
so that the triangle equalities hold:
\[
\begin{tikzcd}[column sep=small] & \cB \arrow[d, phantom, "\scriptstyle\Downarrow\epsilon"] \arrow[dr,  "F" description] \arrow[rr, equals] & \arrow[d, phantom, "\scriptstyle\Downarrow\eta"] & \cB \arrow[dr, phantom, "="] & [+8pt] \cB \arrow[d, phantom, "="] &  & \cB \arrow[dr, "F"'] \arrow[rr, equals] &\arrow[d, phantom, "\scriptstyle\Downarrow\eta"] & \cB \arrow[d, phantom, "\scriptstyle\Downarrow\epsilon"] \arrow[dr, "F"] & \arrow[dr, phantom, "="]  &[+8pt] \cB \arrow[d, phantom, "="] \arrow[d, bend left, start anchor=-65, end anchor=68, "F"] \arrow[d, bend right, start anchor=-115, end anchor=112, "F"'] \\ 
\cA \arrow[ur, "U"] \arrow[rr, equals] & ~& \cA \arrow[ur, "U"']  & ~ & \cA \arrow[u, bend left, start anchor=115, end anchor=-115, "U"] \arrow[u, bend right, start anchor=65,end anchor=-65, "U"']  & & &  \cA \arrow[rr, equals] \arrow[ur,  "U" description] &~ & \cA &  \cA .
\end{tikzcd}
\]
\end{defn}

Using the 2-terminal object $1$, Definition \ref{defn:adjunction} specializes to give notions of an {initial object} or a {terminal object} inside an $\infty$-category $\cA$.

\begin{defn}\label{defn:terminal} An \textbf{initial object} in $\cA$ is given by an adjunction as below-left, while a \textbf{terminal object} in $\cA$ is given by an adjunction as below-right:
\[
\begin{tikzcd}
1 \arrow[r, bend left, "i"] \arrow[r, phantom, "\bot"] & \cA \arrow[l, bend left, "!"] & 1 \arrow[r, bend right, "t"'] \arrow[r, phantom, "\bot"] & \cA . \arrow[l, bend right, "!"']
\end{tikzcd}
\]
\end{defn}

More general limits are defined using the notion of an {absolute right lifting diagram}, which we now recall. Colimits are defined co-dually using absolute left lifting diagrams.

\begin{defn}\label{defn:lifting} In a 2-category, an \textbf{absolute right lifting} of a 1-cell $G \colon \cC \to \cA$ through a 1-cell $F \colon \cB \to \cA$ is given by a 1-cell $R \colon \cC \to \cB$ and a 2-cell
\[
\begin{tikzcd}
 \arrow[dr, phantom, "\scriptstyle\Downarrow\rho" pos=.85] & \cB \arrow[d, "F"] \\ \cC \arrow[r, "G"'] \arrow[ur, "R"] & \cA
 \end{tikzcd}
 \]
 so that any 2-cell as below-left factors uniquely as below-right:
 \begin{equation}\label{eq:abs-lifting-bij}  
\begin{tikzcd}
\cX \arrow[dr, phantom, "\scriptstyle\Downarrow\chi"]\arrow[r, "B"] \arrow[d, "C"'] & \cB \arrow[d, "F"] & \arrow[d, phantom, "="] &  \cX \arrow[dr, phantom, "\scriptstyle\exists!\Downarrow\zeta" pos=.1,  "\scriptstyle\Downarrow\rho" pos=.85 ]\arrow[r, "B"] \arrow[d, "C"'] & \cB \arrow[d, "G"]  \\  \cC \arrow[r, "G"'] & \cA &~ & \cC \arrow[ur, "R" description, pos=.4] \arrow[r, "G"'] & \cA .
\end{tikzcd}
\end{equation}
\end{defn}

The adjective ``absolute'' refers to the following stability property, whose proof is left to the reader:

\begin{lem}\label{lem:absolute-lifting}
 If 
 \[
\begin{tikzcd}
 \arrow[dr, phantom, "\scriptstyle\Downarrow\rho" pos=.85] & \cB \arrow[d, "F"] \\ \cC \arrow[r, "G"'] \arrow[ur, "R"] & \cA
 \end{tikzcd}
 \]
  is an absolute right lifting of $G$ through $F$, then for any $C\colon \cX \to \cC$, $(RC,\rho C)$ is an absolute right lifting of $GC$ through $F$.  \qed
\end{lem} 

For example, by a pasting diagram chase, which is left to the reader:

\begin{lem}\label{lem:counit-as-lifting}
A 2-cell $\epsilon \colon FU \Rightarrow \id_\cA$  defines an absolute right lifting
\[
\begin{tikzcd}
 \arrow[dr, phantom, "\scriptstyle\Downarrow\epsilon" pos=.85] & \cB \arrow[d, "F"] \\ \cA \arrow[r, equals] \arrow[ur, "U"] & \cA
 \end{tikzcd}
 \]
 if and only if it defines the counit of an adjunction $F \dashv U$. \qed
 \end{lem}
 
 For an $\infty$-category $\cA$ in an $\infty$-cosmos and an ordinary category $\cJ$, the simplicial cotensor with the nerve of $\cJ$ defines the $\infty$-category $\cA^\cJ$ of $\cJ$-shaped diagrams in $\cA$. In a cartesian closed $\infty$-cosmos, the exponentials $\cA^\cJ$ also define an $\infty$-category of diagrams, where in this case $\cJ$ is another $\infty$-category. In either context, restriction along the unique map $! \colon \cJ \to 1$ defines a constant diagram functor $\Delta \colon \cA \to \cA^\cJ$ and we say that $\cA$ \textbf{admits all limits of shape $\cJ$} just when this functor admits a right adjoint:
 \begin{equation}\label{eq:all-limits}
 \begin{tikzcd}
 \cA^\cJ \arrow[r, bend right, "\lim"' pos=.45] \arrow[r, phantom, "\bot"] & \cA . \arrow[l, bend right, "\Delta"' pos=.55]
 \end{tikzcd}
 \end{equation}
 The counit $\epsilon$ of such an adjunction defines the data of a limit cone, which restricts to define an absolute right lifting of an individual diagram $D \colon 1 \to \cA^\cJ$
\[
\begin{tikzcd}
&  \arrow[dr, phantom, "\scriptstyle\Downarrow\epsilon" pos=.9] & \cA \arrow[d, "\Delta"] \\ 1 \arrow[r, "D"]  & \cA^\cJ \arrow[ur, "\lim"] \arrow[r, equals]& \cA^\cJ.
 \end{tikzcd}
 \]
 This motivates the following definition.

\begin{defn}\label{defn:limit} A \textbf{limit} of a $\cJ$-shaped diagram $D \colon 1\to \cA^\cJ$ is an absolute right lifting of $d$ through the constant diagram functor $\Delta \colon \cA \to \cA^\cJ$.
\[
\begin{tikzcd}
 \arrow[dr, phantom, "\scriptstyle\Downarrow\lambda" pos=.9] & \cA \arrow[d, "\Delta"] \\ 1 \arrow[r, "D"'] \arrow[ur, "\ell"] & \cA^\cJ.
 \end{tikzcd}
 \]
 Here $\ell \in \cA$ defines the limit object, while $\lambda \colon \Delta \ell \To D$ is the limit cone.
 \end{defn}
 
It follows from a straightforward pasting diagram chase that right adjoints between $\infty$-categories preserve limits in $\infty$-categories \cite[2.4.2]{RV-elements}.

\section[2-categories, 2-functors, and 2-natural transformations]{The formal theory of 2-categories, 2-functors, and 2-natural transformations}\label{sec:2-cat-cosmos}

When interpreted in the 2-category $\twoCat$,
Definitions 3.1 and 3.2 spell out familiar notions of equivalences and adjunctions between 2-categories.

\begin{prop}\label{prop:strict-equivalence} 
The following are equivalent and define what it means for a  2-functor $F \colon \cB \to \cA$ between 2-categories to be a \textbf{2-equivalence}:
\begin{enumerate}
\item $F \colon \cB \to \cA$ is an equivalence in $\twoCat$.
\item $F \colon \cB \to \cA$ satisfies the following pair of properties:
\begin{enumerate}
\item $F$ is surjective on objects up to equivalence, and
\item the induced functor $F \colon \cB(x,y) \to \cA(Fx,Fy)$ is an isomorphism of categories for all $x,y \in \cB$. \qed
\end{enumerate}
\end{enumerate}
\end{prop}

\begin{prop}\label{prop:strict-adjunction}
The following are equivalent and define a \textbf{2-adjunction}:
\begin{enumerate}
\item A pair of 2-functors $F \colon \cB \to \cA$ and $U \colon \cA \to \cB$ between 2-categories form an adjunction $F \dashv U$ in $\twoCat$ with unit $\eta$ and counit $\epsilon$.
\item The 2-natural transformations defined by composition with $\eta$ and $\epsilon$
\[ \begin{tikzcd}\cA(Fb,a) \arrow[r, bend left, "U(-)\cdot \eta"] \arrow[r, phantom, "\cong"] & \cB(b,Ua) \arrow[l, bend left, "\epsilon\cdot F(-)"]
\end{tikzcd}
\]
define inverse isomorphisms of hom-categories for all $a \in \cA$ and $b \in \cB$.
\end{enumerate}
\end{prop}

We seek analogous characterizations of the terminal objects, absolute right lifting diagrams, and limits in $\twoCat$.

\subsection{Terminal objects}

We now explain what it means for a 2-category to admit a terminal object in the sense of Definition \ref{defn:terminal}.
\begin{prop} The following are equivalent and define what it means for a 2-category $\cA$ to admit a \textbf{2-terminal object} $t$:
\begin{enumerate}
\item The unique 2-functor $! \colon \cA \to 1$ admits a right adjoint $t \colon 1 \to \cA$ in $\twoCat$.
\item For each $a \in \cA$, there is a unique 1-cell $a \to t$ admitting only identity endomorphisms. 
\end{enumerate}
\end{prop}
\begin{proof}
By Proposition \ref{prop:strict-adjunction} the adjunctions in $\twoCat$ are exactly the 2-adjunctions. There is a 2-adjunction $! \dashv t$ if and only if there is an object $t \in \cA$ equipped with a 
2-natural isomorphism of hom-categories
\[ \catone \cong \cA(a,t)\]
for all $a \in \cA$. Thus we see that terminal objects are precisely the 2-terminal objects from 2-category theory.
\end{proof}

\subsection{Absolute right lifting diagrams}

By Lemma \ref{lem:counit-as-lifting}, the counit of any adjunction $F \dashv U$ in $\twoCat$ defines an absolute right lifting diagram in $\twoCat$. In fact, this absolute right lifting diagram necessarily satisfies a stronger universal property that is expressible because $\twoCat$ underlies a strict 3-category whose 3-cells are the modifications.

\begin{defn}
  For 2-categories $\cA$, $\cB$, and $\cC$; 2-functors $F \colon \cB \to \cA$, $G \colon \cC \to \cA$, and $R \colon \cC \to \cB$; and a 2-natural transformation $\rho \colon FR \Rightarrow G$ as below
\[
\begin{tikzcd} \arrow[dr, phantom, "\scriptstyle\Downarrow\rho" pos=.85] & \cB \arrow[d, "F"] \\ \cC \arrow[r, "G"'] \arrow[ur, "R"] & \cA
\end{tikzcd}
\]
we say that $(R,\rho)$ defines an \textbf{enriched absolute right lifting diagram} if for all 2-categories $\cX$ and 2-functors $B \colon \cX \to \cB$ and $C \colon \cX \to \cC$ the mapping 
\[
\begin{tikzcd} \cB^{\cX}(B,RC) \arrow[r, "\cong"', "\rho \cdot F(-)"] & \cA^{\cX}(FB,GC)
\end{tikzcd}
\]
defines an isomorphism of categories. 

This means that
\begin{enumerate}
    \item there is a bijection between 2-natural transformations $\zeta \colon B \To RC$ and $\chi \colon FB \To GC$ implemented by pasting with $\rho$ as in \eqref{eq:abs-lifting-bij}, and 
    \item the above bijection extends to a bijection between modifications $\Phi \colon \zeta \Rrightarrow \zeta'$ and $\Psi \colon \chi \Rrightarrow \chi'$.
          \end{enumerate}
  \end{defn}

  The bijection on objects captures precisely the universal property of absolute right lifting diagrams in $\twoCat$, so enriched absolute right liftings are in particular absolute right liftings in $\twoCat$. In the present context, Lemma \ref{lem:counit-as-lifting} extends as follows:

\begin{lem}\label{lem:enriched-counit-as-lifting}
A 2-cell $\epsilon\colon FU \Rightarrow \id_\cA$ defines the counit of a 2-adjunction iff the pair $(U, \epsilon)$ defines an enriched absolute right lifting diagram.
\end{lem}
\begin{proof}
For any 2-category $\cX$, the 2-functor $(-)^\cX \colon \twoCat \to \twoCat$ preserves adjunctions, carrying a 2-adjunction $F \dashv U$ to the 2-adjunction whose corresponding 2-functors are defined by post-composition and post-whiskering: 
\[
\begin{tikzcd}[column sep=large]
\cB^\cX \ar[r, phantom, "\bot"] \ar[r,bend left,"F_*"] & \cA^\cX.\ar[l,bend left,"U_*"]
\end{tikzcd}
\]
The counit $\epsilon_* \colon F_*U_* \To \id_{\cA^\cX}$ of this 2-adjunction is similarly defined by post-composition. For each 2-functor $A \colon \cX \to \cA$, the component of $\epsilon_*$ at $A$ is $\epsilon{A}$. 

By Proposition \ref{prop:strict-adjunction}, this second 2-adjunction induces an isomorphism of categories
\begin{equation}\label{eq:enriched-abs-right-lifting}
\begin{tikzcd} \cB^\cX(B,UA) \arrow[r, "\cong"', "\epsilon \cdot F(-)"] & \cA^\cX(FB,A)
\end{tikzcd}
\end{equation}
 for any 2-functors $A \colon \cX \to \cA$ and $B \colon \cX \to \cB$. This isomorphism is the universal property of the enriched absolute right lifting diagram $(U, \epsilon)$.

Conversely, if $(U, \epsilon)$ defines an enriched absolute right lifting, then the 2-natural map $\epsilon \cdot F(-) \colon \cB^\cX(B,UA) \to \cA^\cX(FB,A)$ defines an isomorphism of categories for $A \in \cA^\cX$ and $B \in \cB^\cX$. When $\cX=1$, this property asserts that $F$ is left 2-adjoint to $U$.
\end{proof}

However, a 2-natural transformation may define an absolute right lifting diagram without also defining an enriched absolute right lifting diagram:

\begin{ex} For any 2-category $\cA$ and $a \in \cA$, the identity 2-cell
\[
\begin{tikzcd} \arrow[dr, phantom, "\scriptstyle\Downarrow\id_a" very near end] & \cA_{0,1} \arrow[d, hook, "I"] \\ 1 \arrow[ur, "a"] \arrow[r, "a"'] & \cA
\end{tikzcd}
\]
defines an absolute right lifting of $a \colon 1 \to \cA$ through the inclusion of the full sub 1-category $I \colon \cA_{0,1} \hookrightarrow \cA$ that contains all of the objects and 1-cells but only identity 2-cells. However, if $\cA$ contains a non-identity 2-cell
\[ \begin{tikzcd} x \arrow[r, bend left, "f"] \arrow[r, bend right, "g"'] \arrow[r, phantom, "\scriptstyle\Downarrow\alpha"] & a
\end{tikzcd}
\]
with codomain $a$, then $\alpha$ defines an arrow in the hom-category $\cA(Ix,a)$ that does not lift to $\cA_{0,1}(x,a)$. Thus this absolute right lifting diagram in $\twoCat$ is not an enriched absolute right lifting diagram.
\end{ex}

\subsection{Limits}

Since $\twoCat$ is cartesian closed, a 2-functor $D \colon \cJ \to \cA$ defines a map $D \colon \catone \to \cA^\cJ$. Inspired by the results of the previous sections, our aim in this section is to compare the universal property expressed by an absolute right lifting diagram 
\[
\begin{tikzcd} \arrow[dr, phantom, "\scriptstyle\Downarrow\lambda" pos=.9] & \cA \arrow[d, "\Delta"] \\ \catone \arrow[r, "D"'] \arrow[ur, "\ell"] & \cA^\cJ
\end{tikzcd}
\]
in $\twoCat$ with the notion of a 2-limit.

\begin{defn}\label{defn:2-limit} A pair $(\ell, \lambda)$ defines a \textbf{2-limit} of the diagram $D \colon \cJ \to \cA$ just when the 2-natural transformation
\[
\lambda \cdot \Delta(-) \colon \cA(a,\ell) \to \cA^\cJ(\Delta a, D)
\]
defines an isomorphism of hom-categories for every $a \in \cA$.
\end{defn}

\begin{ex}\label{ex:all-limit-is-2-limit} 
If $\cA$ admits all $\cJ$-shaped limits, in the sense of the 2-adjunction $\Delta \dashv \lim$ of \eqref{eq:all-limits}, then for any $D \in \cA^\cJ$ the component $(\lim D, \epsilon_D)$ defines a 2-limit with the 2-natural isomorphism of Definition \ref{defn:2-limit} given by adjoint transposition.
\end{ex}

By Lemma \ref{lem:enriched-counit-as-lifting}, if $\cA$ admits all $\cJ$-shaped limits, then as observed in Example \ref{ex:all-limit-is-2-limit} these are 2-limits. They define enriched absolute right lifting diagrams and thus satisfy the $\infty$-categorical criterion of Definition \ref{defn:limit}. In fact, the general notion of 2-limit is precisely captured by enriched absolute right lifting diagrams.

\begin{prop}\label{prop:2-limit-as-enriched-lifting} A 2-natural transformation $\lambda \colon \Delta \ell \To D$ defines a 2-limit if and only if $(\ell,\lambda)$ defines an enriched absolute right lifting of $D$ through $\Delta$.
\end{prop}
\begin{proof}
When $\cX = 1$, the universal property of an enriched absolute right lifting $(\ell,\lambda)$ specializes to the universal property of a 2-limit $\lambda \colon \Delta \ell \To D$. It remains  to show that the universal property of a 2-limit defines an enriched absolute right lifting. To that end, we first establish a bijection between 2-cells of the form
\[
\begin{tikzcd}
\cX \arrow[dr, phantom, "\scriptstyle\Downarrow\chi"]\arrow[r, "A"] \arrow[d, "!"'] & \cA \arrow[d, "\Delta"] & \arrow[d, phantom, "="] &  \cX \arrow[dr, phantom, "\scriptstyle\exists!\Downarrow\zeta" pos=.1,  "\scriptstyle\Downarrow\lambda" pos=.9 ]\arrow[r, "A"] \arrow[d, "!"'] & \cA \arrow[d, "\Delta"]  \\  1 \arrow[r, "D"'] & \cA^\cJ &~ & 1 \arrow[ur, "\ell" description, pos=.4] \arrow[r, "D"'] & \cA^\cJ.
\end{tikzcd}
\]
Given a 2-natural transformation $\chi \colon \Delta A \To D!$ and $x \in \cX$ we define a 1-cell $\zeta_x \colon Ax \to \ell$ in $\cA$ to be the unique 1-cell corresponding to $\chi_x \colon \Delta Ax \to D$ under the isomorphism
\begin{equation}\label{eq:2-limit-UP}
\begin{tikzcd} \cA(Ax,\ell) \arrow[r, "\cong"', "\lambda \cdot \Delta(-)"] & \cA^\cJ(\Delta Ax,D).
\end{tikzcd}
\end{equation}
Since $\chi$ and the isomorphism \eqref{eq:2-limit-UP} are 2-natural, these components $\zeta_x$ assemble into a 2-natural transformation $\zeta \colon A \To \ell!$. Since a 2-natural transformation is determined by its components, $\zeta$ is clearly the unique factorization of $\chi$ through $\lambda$. 

We now use the fact that \eqref{eq:2-limit-UP} is a 2-natural isomorphism of categories, rather than a natural isomorphism of sets, to extend this bijective correspondence to modifications $\Phi \colon \zeta \Rrightarrow \zeta'$ and $\Psi \colon \chi \Rrightarrow \chi'$ where $\chi = \lambda! \cdot \Delta \zeta$ and $\chi' = \lambda! \cdot \Delta \zeta'$. The component of a modification  $\Psi \colon \chi \Rrightarrow \chi'$ at $x \in \cX$ is given by a 2-cell $\Psi_x \colon \chi_x \To \chi'_x$ in $\cA^\cJ$, defining an arrow in the right-hand hom-category of \eqref{eq:2-limit-UP}. Thus we define $\Phi_x \colon \zeta_x \To \zeta'_x$ to be the corresponding arrow in the left-hand hom-category. The final task is to show that these 2-cells assemble into a modification $\Phi \colon \zeta \Rrightarrow \zeta'$, which is again a consequence of the 2-naturality  of \eqref{eq:2-limit-UP} in $x \in \cX$ and the fact that $\Psi$ defines a modification.
\end{proof}

However, if $(\ell,\lambda)$ defines a mere absolute right lifting diagram in $\twoCat$ rather than an \emph{enriched} absolute right lifting diagram, this data defines an unenriched limit but not necessarily an enriched limit, as we now illustrate.

\begin{ex}
Consider the 2-category $\cA$ depicted below:
\[ \begin{tikzcd} a & p \arrow[l, bend left, "\ell"] \arrow[l, bend right, "\ell"'] \arrow[l, phantom, "\scriptstyle\Downarrow\alpha"] \arrow[r, bend left, "r"] \arrow[r, bend right, "r"'] \arrow[r, phantom, "\scriptstyle\Downarrow\beta"] & b
\end{tikzcd}
\]
with three objects and two non-identity arrows, each of which admits a non-identity endomorphism. The object $p$ defines the product of $a$ and $b$ but this product is not a 2-product, since the 1-cell $\id_p$ does not admit any endomorphisms.

We claim, however, that
\[
\begin{tikzcd} \arrow[dr, phantom, "{\scriptstyle\Downarrow(\ell,r)}" pos=.9] & \cA \arrow[d, "\Delta"] \\ 1 \arrow[r, "{(a,b)}"'] \arrow[ur, "p"] & \cA \times \cA
\end{tikzcd}
\]
is an absolute right lifting diagram in $\twoCat$. To see this, first observe that there is a unique functor $A \colon \cX \to \cA$ admitting a 2-natural transformation
\[
\begin{tikzcd} \cX \arrow[r, "A"] \arrow[d, "!"'] \arrow[dr, phantom, "\scriptstyle\Downarrow\chi"] & \cA \arrow[d, "\Delta"] \\ 1 \arrow[r, "{(a,b)}"'] & \cA \times \cA .
\end{tikzcd}
\]
The component of $\chi$ at $x \in \cX$ is given by a pair of 1-cells in $\cA$ from $Ax$ to $a$ and to $b$, and such 1-cells only exist if $Ax = p$. Since $p$ admits no endomorphic 1-cells or 2-cells, this tells us that the 2-functor $A$ is necessarily constant at $p$, in which case the only 2-natural transformation $\chi$ is constant at the pair of 1-cells $(\ell,r)$. So the universal property of the absolute right lifting diagram need only be checked in a single case, for which we have the unique factorization:
\[
\begin{tikzcd}
\cX \arrow[dr, phantom, "\scriptstyle\Downarrow{(\ell,r)}"]\arrow[r, "x!"] \arrow[d, "!"'] & \cA \arrow[d, "\Delta"] & \arrow[d, phantom, "="] &  \cX \arrow[dr, phantom, "\scriptstyle\exists!\Downarrow\id" pos=.1,  "\scriptstyle\Downarrow{(\ell,r)}" pos=.9 ]\arrow[r, "x!"] \arrow[d, "!"'] & \cA \arrow[d, "\Delta"]  \\  1 \arrow[r, "{(a,b)}"'] & \cA \times \cA &~ & 1 \arrow[ur, "p" description, pos=.4] \arrow[r, "{(a,b)}"'] & \cA \times \cA . 
\end{tikzcd}
\]
\end{ex}

There is another hint that Definition \ref{defn:limit} does not specialize to the notion of 2-limit in $\twoCat$. In any $\infty$-cosmos, a limit $(\ell,\lambda)$ of a diagram $D \colon 1 \to \cA^\cJ$ may be characterized as a terminal object in the $\infty$-category $\Hom_{\cA^\cJ}(\Delta,D)$ defined by the pullback
\[
\begin{tikzcd} \Hom_{\cA^\cJ}(\Delta,D) \arrow[r] \arrow[d, two heads] \arrow[dr, phantom, "\lrcorner" very near start] & (\cA^\cJ)^\cattwo \arrow[d, two heads, "{(\cod, \dom)}"] \\ 1 \times \cA \arrow[r, "D \times \Delta"] & \cA^\cJ \times \cA^\cJ.
\end{tikzcd}
\]
In the $\infty$-cosmos $\twoCat$, this is the ``strict slice 2-category'' considered by clingman and Moser \cite[2.7]{clingmanMoser}. While they show that any 2-limit defines a 2-terminal object in the 2-category $\Hom_{\cA^\cJ}(\Delta,D)$, they prove that 2-terminal objects need not necessarily define 2-limits \cite[2.11-12]{clingmanMoser}.

\section[bicategories, normal pseudofunctors, and icons]{The formal theory of bicategories, normal pseudofunctors, and icons}\label{sec:icon-cosmos}

We now turn our attention to the $\infty$-cosmos $\Icon$ of bicategories, normal pseudofunctors, and icons. The equivalences  in the sense of Definition \ref{defn:equivalence} are characterized by an observation of Lack in \cite{Lack-icons}.

\begin{prop}\label{prop:icon-equivalence} 
For a normal pseudofunctor $F \colon \cB \to \cA$ between bicategories, the following are equivalent:
\begin{enumerate}
\item $F \colon \cB \to \cA$ is an equivalence in $\Icon$.
\item $F \colon \cB \to \cA$ satisfies the following pair of properties:
\begin{enumerate}
\item $F$ is bijective on objects, and
\item the induced functor $F \colon \cB(x,y) \to \cA(Fx,Fy)$ is an equivalence of categories for all $x,y \in \cB$. \qed
\end{enumerate}
\end{enumerate}
\end{prop}

In this section, we seek analogous characterizations of the adjunctions, terminal objects, and limits in this 2-category.

\subsection{Adjunctions}

Suppose $F \colon \cB \to \cA$ is left adjoint to $U \colon \cA \to \cB$ in $\Icon$. Then by the presence of the unit and counit icons $\eta\colon \id_\cB \To UF$ and $\epsilon \colon FU \To \id_\cA$, the maps $F$ and $U$ must necessarily define an inverse isomorphism between $\ob \cA$ and $\ob \cB$. Now the hom-category actions of $F$ and $U$ define functors 
\[
\begin{tikzcd} \cB(b,Ua) \arrow[r, bend left, "F"] & \cA(Fb,a) \arrow[l, bend left, "U"] 
\end{tikzcd}
\]
for each $b \in \cB$ and $a \in \cA$, giving rise to the following characterization of adjunctions in $\Icon$.

\begin{prop}\label{prop:icon-adjunction} The following are equivalent:
\begin{enumerate}
\item An adjunction $F \colon \cB \to \cA$ and $U \colon \cA \to \cB$ in $\Icon$ with $F \dashv U$.
\item An inverse  isomorphism $F \colon \ob \cB \cong \ob \cA$, $U \colon \ob \cA \cong \ob \cB$ together with a local adjunction
\[
\begin{tikzcd} \cB(b,Ua) \arrow[r, bend left, "F"] \arrow[r, phantom, "\bot"] & \cA(Fb,a) \arrow[l, bend left, "U"] 
\end{tikzcd}
\]
for all $a \in \cA$ and $b \in \cB$.
\end{enumerate}
\end{prop}

We prove this result as a special case of a more general result appearing in Proposition \ref{prop:icon-lifting}; see Observation \ref{obs:icon-adjunction}.

\subsection{Terminal objects}

We now interpret Definition \ref{defn:terminal} in $\Icon$.

\begin{prop}\label{prop:icon-terminal} A bicategory $\cA$ admits a terminal object $t \colon 1 \to A$ defining a right adjoint to the unique normal pseudofunctor $! \colon A \to 1$ if and only if:
\begin{enumerate}
\item The bicategory $\cA$ has a single object $t$.
\item The identity 1-cell $\id_t \in \cA(t,t)$ is terminal in the unique hom-category of $\cA$.
\end{enumerate}
\end{prop}

Recall a 1-object bicategory is equally a monoidal category. Proposition \ref{prop:icon-terminal} says that a monoidal category, regarded as a 1-object bicategory, admits a terminal object in $\Icon$ if and only if it is \textbf{semi-cartesian}, meaning its monoidal unit is a terminal object.

\begin{proof}
Unpacking Definition \ref{defn:terminal}, a terminal object in a bicategory $\cA$ in $\Icon$ is given by an adjunction
\[
\begin{tikzcd} \cA \arrow[r, bend left, "!"] \arrow[r, phantom, "\bot"] & 1 \arrow[l, bend left, "t"]
\end{tikzcd}
\]
involving the terminal 2-category $1$. By Proposition \ref{prop:icon-adjunction}, this tells us in particular that $\cA$ must be a 1-object bicategory, so $t \colon 1 \to \cA$ is the unique object in $\cA$.

The remaining data defines  a local adjunction
\[
\begin{tikzcd}\cA(t,t) \arrow[r, bend left, "!"] \arrow[r, phantom, "\bot"] & 1(*,*) \arrow[l, bend left, "t"]
\end{tikzcd}
\]
 involving the unique hom-category of $\cA$ and the terminal category $1(*,*) \cong \catone$. Since $t \colon 1 \to \cA$ is a normal pseudofunctor, the local right adjoint $t \colon \catone \to \cA(t,t)$ picks out the identity 1-cell on $t$. The adjointness tells us that the 1-cell  $\id_t \in \cA(t,t)$ is terminal, i.e., for all $f \colon t \to t$ there exists a unique 2-cell $\eta_f \colon f \To \id_t$ in $\cA$. 
\end{proof}

\subsection{Absolute right lifting diagrams}

We now consider an absolute right lifting diagram in $\Icon$:
\[
\begin{tikzcd}
 \arrow[dr, phantom, "\scriptstyle\Downarrow\rho" pos=.85] & \cB \arrow[d, "F"] \\ \cC \arrow[r, "G"'] \arrow[ur, "R"] & \cA .
 \end{tikzcd}
 \]

\begin{prop}\label{prop:icon-lifting} The following are equivalent:
\begin{enumerate}
\item\label{itm:icon-lifting-defn} The pair $(R,\rho)$ defines absolute right lifting of $G \colon \cC \to \cA$ through $F \colon \cB \to \cA$ in $\Icon$.
\item\label{itm:icon-lifting-local} The normal pseudofunctors $R$, $G$, and $F$ give rise to a pullback diagram of object functions
 \[
 \begin{tikzcd} \ob \cC \arrow[r, "R"] \arrow[d, equals] \arrow[dr, phantom, "\lrcorner" very near start] & \ob \cB \arrow[d, "F"] \\ \ob \cC \arrow[r, "G"'] & \ob \cA
 \end{tikzcd}
 \]
and for all $x, y \in \cC$, the induced diagram of hom-categories
\begin{equation}
\begin{tikzcd}\label{eq:induced-diagram}
\arrow[dr, phantom, "\scriptstyle\Downarrow\rho_{x,y}" pos=.85] & \cB(Rx,Ry) \arrow[d, "F"] \\ \cC(x,y) \arrow[r, "G"'] \arrow[ur, "R"]& \cA(Gx,Gy) = \cA(FRx, FRy)
\end{tikzcd}
\end{equation}
is an absolute right lifting diagram in $\Cat$.
\end{enumerate}
\end{prop}

Our proof of Proposition \ref{prop:icon-lifting} makes use of the following lemma.

\begin{lem}\label{lem:icon-2-adjunction}
  There is a 2-adjunction
  \[ \begin{tikzcd} {}^{\catone+\catone/}\Icon \arrow[r, bend right, "\Hom"'] \arrow[r, phantom, "\bot"] & \Cat \arrow[l, bend right, "\Sigma"' pos=.55]
\end{tikzcd}
\]
between the 2-category of bicategories with a chosen pair of objects, basepoint-preserving normal pseudofunctors, and icons and the 2-category of categories whose right adjoint carries a bipointed bicategory $(\cC,x,y)$ to the hom-category $\cC(x,y)$.
\end{lem}
\begin{proof}
  The left 2-adjoint is the fully faithful 2-functor that carries a category $J$ to the 2-category with two objects $0$ and $1$ and whose hom-categories are given by
  \[\Sigma J(0,1) \coloneqq J, \quad \Sigma J(0,0) \coloneqq \Sigma J(1,1) \coloneqq \catone, \quad \text{and} \quad \Sigma J(1,0) \coloneqq \emptyset.\] The component of the counit at $(\cC,x,y)$ is given by the canonical normal pseudo\-functor $\Sigma \cC(x,y) \to \cC$. We leave the verification of the 2-adjointness of this construction to the reader.
\end{proof}

\begin{proof}[Proof of Proposition \ref{prop:icon-lifting}]
First assume that $(R, \rho)$ defines an absolute right lifting of $G$ through $F$, which implies in particular that the 2-functors $FR$ and $G$ agree on objects.

Observe that there exists an icon $\begin{tikzcd} 1 \arrow[r, bend left, "x"] \arrow[r, bend right, "y"'] \arrow[r, phantom, "\scriptstyle\Downarrow\alpha"] & \cA \end{tikzcd}$ if and only if $x = y$ as objects of $\cA$, in which case there is a unique icon from $x$ to $y$ (the identity icon). The universal property of absolute right lifting diagrams specializes to give a bijection between icons as below-left and as below-right
  \[   
\begin{tikzcd}
1 \arrow[dr, phantom, "\scriptstyle\Downarrow\chi"]\arrow[r, "b"] \arrow[d, "c"'] & \cB \arrow[d, "F"] & \arrow[d, phantom, "\leftrightsquigarrow"] &  1 \arrow[dr, phantom, "\scriptstyle\Downarrow\zeta" pos=.1 ]\arrow[r, "b"] \arrow[d, "c"'] & \cB  \\  \cC \arrow[r, "G"'] & \cA &~ & \cC \arrow[ur, "R"'] & ~
\end{tikzcd}
\]
 implemented by pasting with $\rho$.  The former stand in bijection with pairs of objects $(b,c)$ so that $Fb = Gc$, while the latter stand in bijection with objects of $\cC$. Thus we see that if $(R,\rho)$ is absolute right lifting, then the square
 \[
 \begin{tikzcd} \ob \cC \arrow[r, "R"] \arrow[d, equals] \arrow[dr, phantom, "\lrcorner" very near start] & \ob \cB \arrow[d, "F"] \\ \ob \cC \arrow[r, "G"'] & \ob \cA
 \end{tikzcd}
 \]
 is a pullback.
 
 Now let $x,y \in \cC$. Since the forgetful 2-functor ${}^{\catone+\catone/}\Icon \to \Icon$ is locally fully faithful, the absolute right lifting diagram $(R,\rho)$ lifts to define an absolute right lifting diagram in ${}^{\catone+\catone/}\Icon$ starting from the object $(\cC,x,y)$. Right 2-adjoints preserve absolute right lifting diagrams, so it follows that \eqref{eq:induced-diagram} defines an absolute right lifting diagram in $\Cat$ as claimed.

 For the converse, consider an icon
 \[
 \begin{tikzcd}
\cX \ar[r, "B"] \ar[d, swap, "C"] \ar[dr, phantom, "\scriptstyle\Downarrow\chi"] & \cB \ar[d, "F"] \\
\cC \ar[r, swap, "G"] & \cA .
\end{tikzcd}
\]
The presence of the icon $\chi$ demands that $FB$ and $GC$ agree on objects. It follows from the pullback given in \eqref{itm:icon-lifting-local} that $B$ and $RC$ agree on objects as well, which means it is possible to define an icon between these functors.

We use the remaining condition of \eqref{itm:icon-lifting-local} to specify the data of such an icon. For any $x,y \in \cX$, the given absolute right lifting in $\Cat$ factors the component of the icon $\chi$ at $x,y$ as follows:
\[
\begin{tikzcd}
\cX(x, y) \ar[r, "B"] \ar[d, swap, "C"] \ar[dr, phantom, "\scriptstyle\Downarrow\chi_{x,y}"] & \cB(Bx, By) \ar[d, "F"] \\
\cC(Cx, Cy) \ar[r, swap, "G"] & \cA(FBx,FBy)
\end{tikzcd}
= 
\begin{tikzcd}
  \cX(x, y) \ar[r, "B"] \ar[d, swap, "C"] \ar[dr, phantom, "\scriptstyle\exists!\Downarrow\zeta_{x,y}" very near start] \ar[dr, phantom, "\scriptstyle\Downarrow\rho_{Cx,Cy}" very near end] & \cB(Bx, By) \ar[d, "F"] \\
\cC(Cx, Cy) \ar[r, swap, "G"] \ar[ur, "R" description, pos=.35] & \cA(FBx,FBy) .
\end{tikzcd}
\]
We claim that this defines the component at $x, y \in \cX$ of an icon $\zeta$ in the sense of Definition \ref{defn:icon}. Since icons are determined entirely by their components, it will follow immediately that $\chi = \rho C \cdot F\zeta$ and that this is the unique such factorization, proving that $(R,\rho)$ is absolute right lifting in $\Icon$.

To see that the components $\zeta_{x,y}$ assemble into an icon observe first that since $\chi$ and $\rho$ are icons, for each $x \in \cX$ their components at identity 1-cells define identity natural transformations in $\Cat$. The defining factorization property tells us that $\zeta_{x,x}$ also has its component at $1_{x}$ given by the identity natural transformation. This verifies the unit property.

The final icon condition requires us to verify the pasting equality
\[
  \begin{tikzcd}[sep=tiny]
    \cX(y,z) \times \cX(x,y) \arrow[dr, "\circ"] \arrow[dd, "C \times C"']\arrow[rr, "B \times B"] & & \cB(By,Bz) \times \cB(Bx,By) \arrow[dr, "\circ"] \\
    & \cX(x,z) \arrow[rr, "B"] \arrow[dd, "C"']  \arrow[ur, phantom, "\scriptstyle\cong"] \arrow[dl, phantom, "\scriptstyle\cong"]& & \cB(Bx,Bz) \\
    \cC(Cy, Cz) \times \cC(Cx,Cy) \arrow[dr, "\circ"']& & \arrow[ul, phantom, "\scriptstyle\Downarrow\zeta_{x,z}"] \\
    & \cC(Cx,Cz) \arrow[uurr, "R"'] & & =
    \end{tikzcd}
\]
\[
  \begin{tikzcd}[sep=tiny]
    \cX(y,z) \times \cX(x,y)  \arrow[dd, "C \times C"']\arrow[rr, "B \times B"] & & \cB(By,Bz) \times \cB(Bx,By)  \arrow[dr, "\circ"] \\
    & \arrow[ul, phantom, "{\scriptstyle\Downarrow \zeta_{y,z}\times \zeta_{x,y}}"] & & \cB(Bx,Bz) \\
    \cC(Cy, Cz) \times \cC(Cx,Cy) \arrow[uurr, "R \times R"'] \arrow[dr, "\circ"'] & \arrow[r, phantom, "\scriptstyle\cong"] & ~ \\
    & \cC(Cx,Cz) \arrow[uurr, "R"']
    \end{tikzcd}
\]
for any $x,y,z \in \cX$, where the unlabelled isomorphisms are the coherences of the normal pseudofunctors $B$, $C$, and $R$. By the universal property of the absolute right lifting diagrams in $\Cat$ it suffices to verify this equality after pasting with $\rho_{Cx,Cz}$. At this point, the required pasting equality reduces to the corresponding icon condition for $\chi$.
\end{proof}

\begin{obs}\label{obs:icon-adjunction} Proposition \ref{prop:icon-adjunction} can be deduced as a special case of Lemma \ref{prop:icon-lifting} via Lemma \ref{lem:counit-as-lifting}. An icon $\epsilon \colon FU \To \id_\cA$ defines the counit of an adjunction if and only if it defines an absolute right lifting, which we have just shown is the case if and only if
 \[
 \begin{tikzcd} \ob \cA \arrow[r, "U"] \arrow[d, equals] \arrow[dr, phantom, "\lrcorner" very near start] & \ob \cB \arrow[d, "F"] \\ \ob \cA \arrow[r, equals] & \ob \cA
 \end{tikzcd}
 \]
 is a pullback and for all $z,a \in \cA$, 
 \[
\begin{tikzcd} \arrow[dr, phantom, "\scriptstyle\Downarrow\epsilon" pos=.85] & \cB(Uz,Ua) \arrow[d, "F"] \\ \cA(z,a) \arrow[r, equals] \arrow[ur, "U"]& \cA(z,a)
\end{tikzcd}
\]
is an absolute right lifting diagram of categories. The former property holds if and only if $F \colon \ob \cB \to \ob \cA$ and $U \colon \ob \cA \to \ob \cB$ define an inverse isomorphism. In particular, we may replace $z \in \cA$ by the unique $b \in \cB$ so that $Fb = z$. By Lemma \ref{lem:counit-as-lifting} again, 
 \[
\begin{tikzcd} \arrow[dr, phantom, "\scriptstyle\Downarrow\epsilon" pos=.85] & \cB(b,Ua) \arrow[d, "F"] \\ \cA(Fb,a) \arrow[r, equals] \arrow[ur, "U"]& \cA(Fb,a)
\end{tikzcd}
\]
defines an absolute right lifting if and only if $\epsilon$ is the counit of an adjunction
\[
\begin{tikzcd} \cB(b,Ua) \arrow[r, bend left, "F"] \arrow[r, phantom, "\bot"] & \cA(Fb,a). \arrow[l, bend left, "U"] 
\end{tikzcd}
\]
\end{obs}

\subsection{Limits}

While the 1-category of bicategories and normal pseudofunctors is cartesian closed, the 2-category of bicategories, normal pseudofunctors, and icons is not cartesian closed as a 2-category. A parallel pair of normal pseudofunctors $F,G \colon \cA \to \cB$ may admit a (non-identity) icon between them, but the transposed normal pseudofunctors $F,G\colon 1 \to \cB^{\cA}$ only admit an icon between them when $F=G$, in which case the only icon is the identity. 

However, since the 2-category $\Icon$ has PIE-limits, it necessarily also has cotensors by any 1-category. The cotensor of a bicategory $\cA$ with a 1-category $\cJ$ is  the bicategory $\cA^\cJ$ characterized by the natural isomorphism of categories
\[\Icon(\cX,\cA^{\cJ})\cong \Icon(\cX,\cA)^{\cJ}.\]
By substituting the free 0-cell $1$, free 1-cell $C_1$, and free 2-cell $C_2$ for $\cX$, we can describe the bicategory $\cA^{\cJ}$.

\begin{obs}\label{obs:icon-cotensors}
  For a 1-category $\cJ$ and bicategory $\cA$:
  \begin{itemize}
  \item Objects in $\cA^{\cJ}$ correspond to functors $\cJ \to \Icon(1,\cA)$. Since the hom-category $\Icon(1,\cA)$ is isomorphic to the discrete category on the set of objects of $\cA$, the objects in $\cA^{\cJ}$ can also be understood as functions $\pi_{0}\cJ \to \ob{\cA}$ from the set of path components of $\cJ$ to the set of objects of $\cA$.
  \item 1-cells in $\cA^{\cJ}$ correspond to functors $\cJ \to \Icon(C_1, \cA)$. These in turn can be understood via the isomorphism
    \[ \Icon(C_1,\cA) \cong \coprod\limits_{x,y \in \cA}\cA(x,y).\]
    The source and target of a 1-cell $\phi \colon S \to T$ in $\cA^{\cJ}$ are, of course, objects in $\cA^{\cJ}$, these being given by functions $S,T \colon \pi_{0}\cJ \to \ob{\cA}$. The rest of the data of the 1-cell $\phi$ is then given by a ``dependent functor" that sends each $j \in \cJ$ to a 1-cell in $\cA(Sj, Tj)$ and a morphism $f \colon j \to k \in \cJ$ to a 2-cell between such 1-cells; note that since $j$ and $k$ belong to the same path component, $Sj=Sk$ and $Tj=Tk$.
  \item 2-cells in $\cA^{\cJ}$ correspond to functors $\cJ \to \Icon(C_{2},\cA)$, which can be understood via the isomorphism
    \[\Icon(C_{2}, \cA) \cong \coprod\limits_{x,y \in \cA}\cA(x,y)^{\cattwo}.\]
  \end{itemize}
 \end{obs}

 We use Proposition \ref{prop:icon-lifting} to describe limits in $\Icon$ of 1-category indexed diagrams valued in a bicategory. By Observation \ref{obs:icon-cotensors}, a ``diagram" $D \colon 1 \to \cA^{\cJ}$ is given by a function $D \colon \pi_{0}\cJ \to \ob{\cA}$. Thus the diagrams considered in this context are rather degenerate.

 \begin{prop}\label{prop:icon-limit} A diagram $D\colon 1 \to \cA^{\cJ}$ in $\Icon$ admits a limit if and only if 
 \begin{enumerate}
 \item the diagram $D$ is constant at some object $\ell \in \cA$, and 
 \item $\id_{\ell} \in \cA(\ell,\ell)$ is the limit of the constant $\cJ$-shaped diagram in $\cA(\ell,\ell)$ at this object, with the identity natural transformation defining the limit cone. 
 \end{enumerate} \end{prop}
\begin{proof}
  If $D \colon 1 \to \cA^{\cJ}$ admits a limit, then it must admit a limit cone, this being an icon
  \[
    \begin{tikzcd}
      \arrow[dr, phantom, "\scriptstyle\Downarrow\lambda" very near end] & \cA \arrow[d, "\Delta" ] \\ 1\arrow[ur, "\ell"] \arrow[r, "D"'] & \cA^{\cJ}.
    \end{tikzcd}
  \]
  As observed above, such icons exist if and only if $\Delta \ell = D$. This tells us that $D$ must be the constant diagram at the object $\ell$, in which case the icon $\lambda$ is necessarily the identity icon.

  It remains to describe conditions that make $(\ell, \id_{\Delta\ell})$ an absolute right lifting of the constant diagram $D = \Delta \ell$ through $\Delta$. By Proposition \ref{prop:icon-lifting}, $(\ell,\id_{\Delta\ell})$ defines an absolute right lifting in $\Icon$ if and only if the induced object functions define a pullback diagram --- which is automatic in this case since $\Delta \colon \ob{\cA} \to \ob{\cA^{\cJ}}$ is a monomorphism --- and if the diagram of categories
  \[
    \begin{tikzcd} \arrow[dr, phantom, "\scriptstyle\Downarrow\id" very near end] & \cA(\ell,\ell) \arrow[d] \\ 1 \arrow[r, "\id_{\ell}"'] \arrow[ur, "\id_{\ell}"] & \cA^{\cJ}(\Delta \ell,\Delta \ell) \cong \cA(\ell,\ell)^{\cJ}
      \end{tikzcd}
    \]
    is an absolute right lifting diagram where the isomorphism in the lower-right corner is a consequence of Observation \ref{obs:icon-cotensors}. Thus we see that a constant diagram $D= \Delta \ell$ admits a limit if and only if the identity cone defines a limit for the object $\id_{\ell}$ in the hom-category $\cA(\ell,\ell)$. 
\end{proof}

The conditions of Proposition \ref{prop:icon-limit} are satisfied in particular if either:
   \begin{enumerate}
   \item $\cJ$ has a single path component, or
   \item $\id_{\ell} \in \cA(\ell,\ell)$ is terminal.
   \end{enumerate}
We do not know whether there are any contexts in which the notion of limit captured by this result is of interest.

\section{Another \texorpdfstring{$\infty$}{infinity}-cosmos of bicategories}\label{sec:conjectures}

The most general notion of sameness between 2-categories or bicategories is not captured by either of the $\infty$-cosmoi $\twoCat$ or $\Icon$. A normal pseudofunctor $F \colon \cB \to \cA$ is a \textbf{biequivalence} just when
\begin{enumerate}
\item $F$ is surjective on objects up to isomorphism, and
\item the induced functor $F \colon \cB(x,y) \to \cA(Fx,Fy)$ is an equivalence of categories for all $x,y \in \cB$,
\end{enumerate}
combining the weaker of the two properties captured by the equivalences in Propositions \ref{prop:strict-equivalence} and \ref{prop:icon-equivalence}. In this section, we describe another $\infty$-cosmos of bicategories whose equivalences are the biequivalences due to Alexander Campbell, which he obtains by applying \cite[E.1.3]{RV-elements} to a model structure he constructs on a closely related category of algebraically cofibrant 2-categories \cite{Campbell-talk}.

The underlying 1-category of the $\infty$-cosmos $\Bicat$ is the category of bicategories and normal pseudofunctors. The isofibrations in $\Bicat$ are the \textbf{equifibrations}, normal pseudofunctors $F \colon \cA \to\cB$ with the property that:
\begin{enumerate}
\item Any equivalence $Fa \simeq b$ in $\cB$ can be lifted to an equivalence $a \simeq x$ in $\cA$.
\item For all $x,y \in \cA$, the functor $F \colon \cA(x,y) \to \cB(Fx,Fy)$ is an isofibration of categories.
\end{enumerate}

Unlike the $\infty$-cosmoi considered above, the $\infty$-cosmos $\Bicat$ is not isomorphic to its homotopy 2-category. Here the simplicial functor-space between bicategories $\cA$ and $\cB$ is the the quasi-category whose $n$-simplices are normal pseudofunctors $\cA \times [n] \to \cB$, where $[n]$ is the free 1-category on $n$ composable arrows. Thus, the homotopy 2-category of $\Bicat$ is the 2-category of bicategories, normal pseudofunctors, and equivalence classes of ``enhanced pseudonatural transformations.''\footnote{An \textbf{enhanced pseudonatural transformation} is an a normal pseudofunctor $\cA \times C_1 \to \cB$. Two enhanced pseudonatural transformation represent the same 2-cell in the homotopy 2-category of $\Bicat$ if they bound a normal pseudofunctor $\cA \times C_2 \to \cB$.}
This homotopy 2-category can also be obtained in another fashion using the fact that the category of bicategories and normal pseudofunctors is cartesian closed.  There is a product-preserving functor that takes a bicategory to the category with the same objects and whose arrows are isomorphism classes of parallel 1-cells. Change-of-base then converts the bicategorically enriched category of bicategories to a 2-category of bicategories, namely the homotopy 2-category described above.

Unlike the familiar 2-categories $\twoCat$ and $\Icon$, we are not sure whether this homotopy 2-category of bicategories, normal pseudofunctors, and equivalence classes of enhanced pseudonatural transformations has been studied before. Nevertheless, it would be interesting to explore the formal category theory of the $\infty$-cosmos $\Bicat$. We leave this unfinished end for future work.

\end{document}